\documentclass[11pt]{article}
\usepackage[a4paper,margin=2.6cm]{geometry}
\usepackage{amsmath,amssymb,amsthm}
\usepackage{booktabs}
\usepackage{tikz}
\usetikzlibrary{backgrounds}
\usepackage[hidelinks]{hyperref}

\newtheorem{theorem}{Theorem}
\newtheorem{lemma}[theorem]{Lemma}
\theoremstyle{remark}
\newtheorem{remark}[theorem]{Remark}

\title{A $43$-vertex $4$-chromatic unit-distance graph in the plane with girth $5$}
\author{Chandan Shukla\\ \small\texttt{shukla.chandan12@gmail.com}}
\date{}

\begin{document}
\maketitle

\begin{abstract}
We exhibit a set of $43$ points in the Euclidean plane whose unit-distance graph has $83$ edges, girth $5$ and chromatic number $4$. The smallest previously known $4$-chromatic unit-distance graph of girth $5$, due to Hochberg and O'Donnell (1996), has $45$ vertices. The new graph is $4$-critical, has trivial automorphism group, and has exactly $2n-3$ edges; it is a Laman graph and its unit-distance realisation is infinitesimally rigid. The coordinates are algebraic numbers, specified as the unique real solution of an explicit system of $82$ quadratic equations inside an explicit rational box. Every claim is certified by exact computation: existence and uniqueness of the realisation by a rational contraction-mapping certificate (independently re-checked with interval Krawczyk arithmetic), and non-$3$-colourability both by exhaustive search and by a DRAT proof verified with \texttt{drat-trim}.
\end{abstract}


\section{Introduction}

A \emph{unit-distance graph} is a graph whose vertices are points of the plane, two points being adjacent when their Euclidean distance is exactly~$1$. The Hadwiger--Nelson problem asks for the chromatic number $\chi(\mathbb{R}^2)$ of the graph on all points of the plane; it is known that $5\le\chi(\mathbb{R}^2)\le 7$~\cite{deGrey,Soifer}.

Erd\H{o}s asked whether $4$-chromatic unit-distance graphs can avoid short cycles. Wormald~\cite{Wormald} constructed a $4$-chromatic unit-distance graph of girth~$5$ with $6448$ vertices. O'Donnell later showed that $4$-chromatic unit-distance graphs of every girth exist~\cite{ODonnell}. For small orders, Hochberg and O'Donnell~\cite{HO} found $4$-chromatic unit-distance graphs of girth~$4$ and girth~$5$ with $23$ and $45$ vertices respectively. The girth-$4$ record has since been lowered to $17$ vertices by Exoo and Ismailescu~\cite{EI}; we are not aware of any improvement of the $45$-vertex girth-$5$ example.

Without the geometric restriction, the smallest $4$-chromatic graphs of girth at least~$5$ have $21$ vertices (Royle; see~\cite{EG}). So the minimum order of a $4$-chromatic unit-distance graph of girth~$5$ lies between $21$ and $45$.

In this note we lower the upper bound to $43$.

\begin{theorem}\label{thm:main}
There is a set $P$ of $43$ points in $\mathbb{R}^2$ whose unit-distance graph $G_{43}$ has the following properties:
\begin{enumerate}
\item $G_{43}$ has $43$ vertices and $83$ edges, and its girth is $5$;
\item $\chi(G_{43})=4$;
\item $G_{43}$ is $4$-critical: deleting any vertex or any edge leaves a $3$-colourable graph;
\item $G_{43}$ has degree sequence $3^{31}5^{3}6^{5}7^{4}$ and trivial automorphism group;
\item $G_{43}$ has $2|V|-3$ edges and is a Laman graph, and the realisation $P$ is infinitesimally rigid.
\end{enumerate}
\end{theorem}

The graph was found by a computer search; this note is concerned only with the graph and with the certificate of its properties.

\section{The graph and its realisation}

Label the vertices $0,1,\dots,42$. The $83$ edges of $G_{43}$ are
\begin{center}\small
\begin{minipage}{0.95\textwidth}\raggedright
$0$--$2$, $0$--$28$, $0$--$29$, $1$--$25$, $1$--$30$, $1$--$31$, $2$--$12$, $2$--$31$,
$3$--$6$, $3$--$24$, $3$--$28$, $4$--$5$, $4$--$13$, $4$--$15$, $5$--$7$, $5$--$10$,
$6$--$11$, $6$--$12$, $7$--$11$, $7$--$14$, $8$--$10$, $8$--$13$, $8$--$18$, $9$--$17$,
$9$--$20$, $9$--$21$, $9$--$27$, $9$--$36$, $10$--$16$, $10$--$23$, $10$--$27$, $10$--$28$,
$10$--$31$, $11$--$16$, $11$--$20$, $11$--$29$, $11$--$32$, $12$--$23$, $12$--$33$, $12$--$41$,
$13$--$39$, $14$--$18$, $14$--$39$, $15$--$18$, $15$--$19$, $15$--$31$, $15$--$32$, $15$--$33$,
$16$--$17$, $16$--$25$, $16$--$26$, $16$--$33$, $16$--$35$, $17$--$19$, $19$--$22$, $20$--$31$,
$20$--$37$, $20$--$39$, $21$--$28$, $21$--$35$, $22$--$28$, $22$--$37$, $23$--$34$, $24$--$29$,
$24$--$33$, $25$--$40$, $26$--$30$, $26$--$38$, $27$--$41$, $28$--$32$, $28$--$39$, $30$--$39$,
$31$--$38$, $32$--$34$, $32$--$40$, $32$--$41$, $33$--$36$, $33$--$39$, $34$--$36$, $35$--$37$,
$38$--$42$, $39$--$42$, $40$--$42$.
\end{minipage}
\end{center}

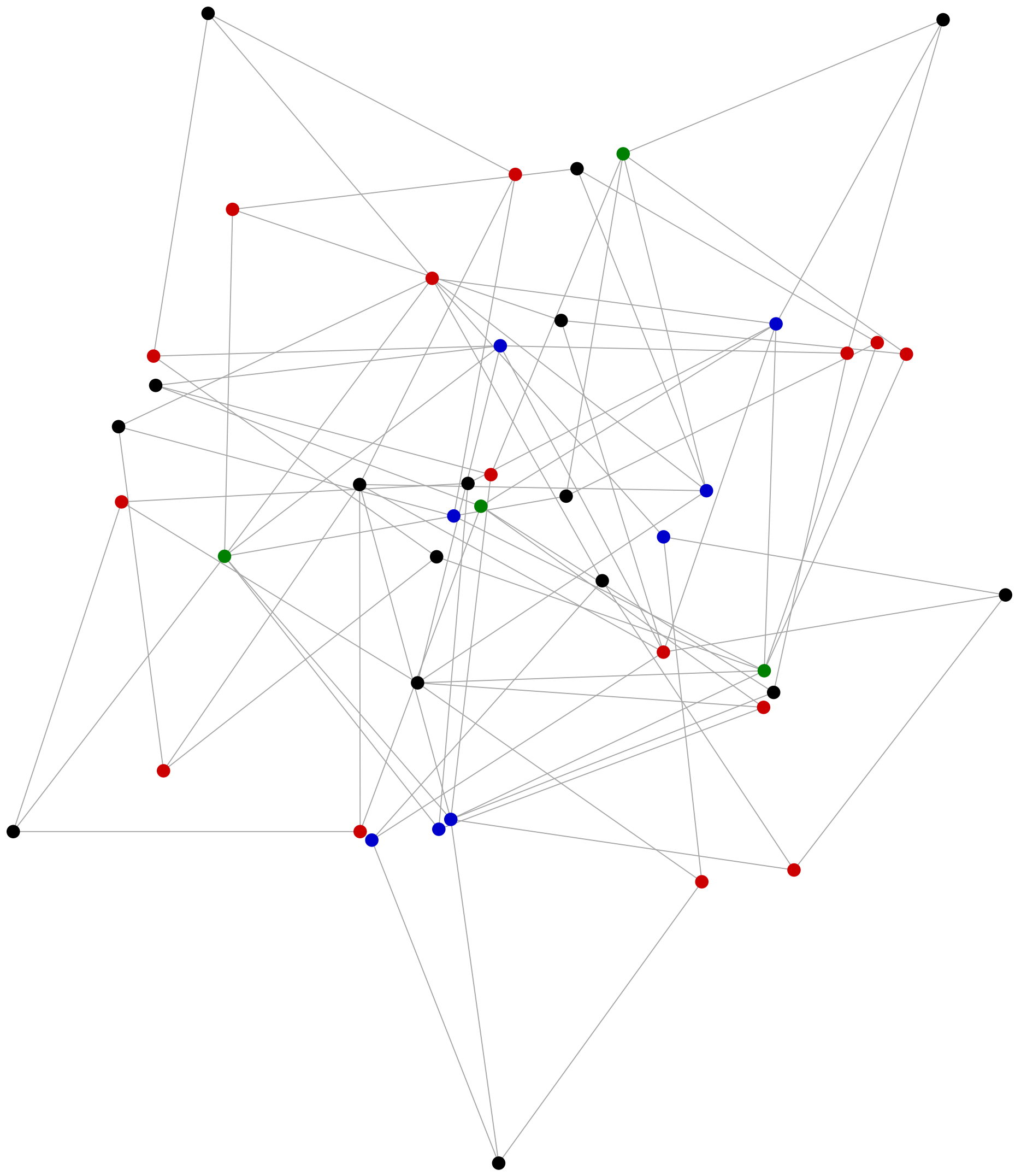
\begin{figure}[t]
\centering
\begin{tikzpicture}[scale=4.2]
\draw[gray!70,line width=0.35pt] (0.00000,0.00000)--(1.00000,0.00000);
\draw[gray!70,line width=0.35pt] (0.00000,0.00000)--(0.60903,0.79314);
\draw[gray!70,line width=0.35pt] (0.00000,0.00000)--(0.31223,0.95001);
\draw[gray!70,line width=0.35pt] (0.43309,0.17518)--(1.22035,0.79180);
\draw[gray!70,line width=0.35pt] (0.43309,0.17518)--(0.30353,1.16675);
\draw[gray!70,line width=0.35pt] (0.43309,0.17518)--(0.99848,1.00000);
\draw[gray!70,line width=0.35pt] (1.00000,0.00000)--(1.34805,0.93748);
\draw[gray!70,line width=0.35pt] (1.00000,0.00000)--(0.99848,1.00000);
\draw[gray!70,line width=0.35pt] (1.22675,0.00674)--(2.16303,0.35799);
\draw[gray!70,line width=0.35pt] (1.22675,0.00674)--(1.31081,1.00320);
\draw[gray!70,line width=0.35pt] (1.22675,0.00674)--(0.60903,0.79314);
\draw[gray!70,line width=0.35pt] (1.03358,-0.02458)--(1.39938,-0.95527);
\draw[gray!70,line width=0.35pt] (1.03358,-0.02458)--(1.69799,0.72279);
\draw[gray!70,line width=0.35pt] (1.03358,-0.02458)--(1.87416,0.51711);
\draw[gray!70,line width=0.35pt] (1.39938,-0.95527)--(1.98493,-0.14463);
\draw[gray!70,line width=0.35pt] (1.39938,-0.95527)--(1.26130,0.03515);
\draw[gray!70,line width=0.35pt] (2.16303,0.35799)--(1.16553,0.42858);
\draw[gray!70,line width=0.35pt] (2.16303,0.35799)--(1.34805,0.93748);
\draw[gray!70,line width=0.35pt] (1.98493,-0.14463)--(1.16553,0.42858);
\draw[gray!70,line width=0.35pt] (1.98493,-0.14463)--(1.87461,0.84926);
\draw[gray!70,line width=0.35pt] (2.25061,-0.11064)--(1.26130,0.03515);
\draw[gray!70,line width=0.35pt] (2.25061,-0.11064)--(1.69799,0.72279);
\draw[gray!70,line width=0.35pt] (2.25061,-0.11064)--(2.86050,0.68185);
\draw[gray!70,line width=0.35pt] (1.75822,1.95291)--(2.57478,1.37564);
\draw[gray!70,line width=0.35pt] (1.75822,1.95291)--(1.99832,0.98216);
\draw[gray!70,line width=0.35pt] (1.75822,1.95291)--(1.59389,0.96650);
\draw[gray!70,line width=0.35pt] (1.75822,1.95291)--(1.37697,1.02844);
\draw[gray!70,line width=0.35pt] (1.75822,1.95291)--(2.68058,2.33925);
\draw[gray!70,line width=0.35pt] (1.26130,0.03515)--(2.16492,0.46349);
\draw[gray!70,line width=0.35pt] (1.26130,0.03515)--(2.19197,0.40101);
\draw[gray!70,line width=0.35pt] (1.26130,0.03515)--(1.37697,1.02844);
\draw[gray!70,line width=0.35pt] (1.26130,0.03515)--(0.60903,0.79314);
\draw[gray!70,line width=0.35pt] (1.26130,0.03515)--(0.99848,1.00000);
\draw[gray!70,line width=0.35pt] (1.16553,0.42858)--(2.16492,0.46349);
\draw[gray!70,line width=0.35pt] (1.16553,0.42858)--(1.99832,0.98216);
\draw[gray!70,line width=0.35pt] (1.16553,0.42858)--(0.31223,0.95001);
\draw[gray!70,line width=0.35pt] (1.16553,0.42858)--(1.40406,1.39972);
\draw[gray!70,line width=0.35pt] (1.34805,0.93748)--(2.19197,0.40101);
\draw[gray!70,line width=0.35pt] (1.34805,0.93748)--(2.19888,1.46292);
\draw[gray!70,line width=0.35pt] (1.34805,0.93748)--(0.41060,1.28559);
\draw[gray!70,line width=0.35pt] (1.69799,0.72279)--(1.20754,1.59426);
\draw[gray!70,line width=0.35pt] (1.87461,0.84926)--(2.86050,0.68185);
\draw[gray!70,line width=0.35pt] (1.87461,0.84926)--(1.20754,1.59426);
\draw[gray!70,line width=0.35pt] (1.87416,0.51711)--(2.86050,0.68185);
\draw[gray!70,line width=0.35pt] (1.87416,0.51711)--(1.57950,1.47271);
\draw[gray!70,line width=0.35pt] (1.87416,0.51711)--(0.99848,1.00000);
\draw[gray!70,line width=0.35pt] (1.87416,0.51711)--(1.40406,1.39972);
\draw[gray!70,line width=0.35pt] (1.87416,0.51711)--(2.19888,1.46292);
\draw[gray!70,line width=0.35pt] (2.16492,0.46349)--(2.57478,1.37564);
\draw[gray!70,line width=0.35pt] (2.16492,0.46349)--(1.22035,0.79180);
\draw[gray!70,line width=0.35pt] (2.16492,0.46349)--(1.26987,0.90946);
\draw[gray!70,line width=0.35pt] (2.16492,0.46349)--(2.19888,1.46292);
\draw[gray!70,line width=0.35pt] (2.16492,0.46349)--(2.49069,1.40894);
\draw[gray!70,line width=0.35pt] (2.57478,1.37564)--(1.57950,1.47271);
\draw[gray!70,line width=0.35pt] (1.57950,1.47271)--(0.63214,1.79288);
\draw[gray!70,line width=0.35pt] (1.99832,0.98216)--(0.99848,1.00000);
\draw[gray!70,line width=0.35pt] (1.99832,0.98216)--(1.62526,1.90997);
\draw[gray!70,line width=0.35pt] (1.99832,0.98216)--(1.20754,1.59426);
\draw[gray!70,line width=0.35pt] (1.59389,0.96650)--(0.60903,0.79314);
\draw[gray!70,line width=0.35pt] (1.59389,0.96650)--(2.49069,1.40894);
\draw[gray!70,line width=0.35pt] (0.63214,1.79288)--(0.60903,0.79314);
\draw[gray!70,line width=0.35pt] (0.63214,1.79288)--(1.62526,1.90997);
\draw[gray!70,line width=0.35pt] (2.19197,0.40101)--(2.40383,1.37831);
\draw[gray!70,line width=0.35pt] (1.31081,1.00320)--(0.31223,0.95001);
\draw[gray!70,line width=0.35pt] (1.31081,1.00320)--(2.19888,1.46292);
\draw[gray!70,line width=0.35pt] (1.22035,0.79180)--(0.40450,1.37007);
\draw[gray!70,line width=0.35pt] (1.26987,0.90946)--(0.30353,1.16675);
\draw[gray!70,line width=0.35pt] (1.26987,0.90946)--(1.44740,1.89357);
\draw[gray!70,line width=0.35pt] (1.37697,1.02844)--(0.41060,1.28559);
\draw[gray!70,line width=0.35pt] (0.60903,0.79314)--(1.40406,1.39972);
\draw[gray!70,line width=0.35pt] (0.60903,0.79314)--(1.20754,1.59426);
\draw[gray!70,line width=0.35pt] (0.30353,1.16675)--(1.20754,1.59426);
\draw[gray!70,line width=0.35pt] (0.99848,1.00000)--(1.44740,1.89357);
\draw[gray!70,line width=0.35pt] (1.40406,1.39972)--(2.40383,1.37831);
\draw[gray!70,line width=0.35pt] (1.40406,1.39972)--(0.40450,1.37007);
\draw[gray!70,line width=0.35pt] (1.40406,1.39972)--(0.41060,1.28559);
\draw[gray!70,line width=0.35pt] (2.19888,1.46292)--(2.68058,2.33925);
\draw[gray!70,line width=0.35pt] (2.19888,1.46292)--(1.20754,1.59426);
\draw[gray!70,line width=0.35pt] (2.40383,1.37831)--(2.68058,2.33925);
\draw[gray!70,line width=0.35pt] (2.49069,1.40894)--(1.62526,1.90997);
\draw[gray!70,line width=0.35pt] (1.44740,1.89357)--(0.56161,2.35765);
\draw[gray!70,line width=0.35pt] (1.20754,1.59426)--(0.56161,2.35765);
\draw[gray!70,line width=0.35pt] (0.40450,1.37007)--(0.56161,2.35765);
\fill[black] (0.00000,0.00000) circle (0.55pt);
\fill[red!80!black] (0.43309,0.17518) circle (0.55pt);
\fill[red!80!black] (1.00000,0.00000) circle (0.55pt);
\fill[blue!80!black] (1.22675,0.00674) circle (0.55pt);
\fill[blue!80!black] (1.03358,-0.02458) circle (0.55pt);
\fill[black] (1.39938,-0.95527) circle (0.55pt);
\fill[red!80!black] (2.16303,0.35799) circle (0.55pt);
\fill[red!80!black] (1.98493,-0.14463) circle (0.55pt);
\fill[red!80!black] (2.25061,-0.11064) circle (0.55pt);
\fill[green!50!black] (1.75822,1.95291) circle (0.55pt);
\fill[blue!80!black] (1.26130,0.03515) circle (0.55pt);
\fill[black] (1.16553,0.42858) circle (0.55pt);
\fill[green!50!black] (1.34805,0.93748) circle (0.55pt);
\fill[black] (1.69799,0.72279) circle (0.55pt);
\fill[blue!80!black] (1.87461,0.84926) circle (0.55pt);
\fill[red!80!black] (1.87416,0.51711) circle (0.55pt);
\fill[green!50!black] (2.16492,0.46349) circle (0.55pt);
\fill[red!80!black] (2.57478,1.37564) circle (0.55pt);
\fill[black] (2.86050,0.68185) circle (0.55pt);
\fill[black] (1.57950,1.47271) circle (0.55pt);
\fill[blue!80!black] (1.99832,0.98216) circle (0.55pt);
\fill[black] (1.59389,0.96650) circle (0.55pt);
\fill[red!80!black] (0.63214,1.79288) circle (0.55pt);
\fill[black] (2.19197,0.40101) circle (0.55pt);
\fill[black] (1.31081,1.00320) circle (0.55pt);
\fill[black] (1.22035,0.79180) circle (0.55pt);
\fill[blue!80!black] (1.26987,0.90946) circle (0.55pt);
\fill[red!80!black] (1.37697,1.02844) circle (0.55pt);
\fill[green!50!black] (0.60903,0.79314) circle (0.55pt);
\fill[red!80!black] (0.31223,0.95001) circle (0.55pt);
\fill[black] (0.30353,1.16675) circle (0.55pt);
\fill[black] (0.99848,1.00000) circle (0.55pt);
\fill[blue!80!black] (1.40406,1.39972) circle (0.55pt);
\fill[blue!80!black] (2.19888,1.46292) circle (0.55pt);
\fill[red!80!black] (2.40383,1.37831) circle (0.55pt);
\fill[red!80!black] (2.49069,1.40894) circle (0.55pt);
\fill[black] (2.68058,2.33925) circle (0.55pt);
\fill[black] (1.62526,1.90997) circle (0.55pt);
\fill[red!80!black] (1.44740,1.89357) circle (0.55pt);
\fill[red!80!black] (1.20754,1.59426) circle (0.55pt);
\fill[red!80!black] (0.40450,1.37007) circle (0.55pt);
\fill[black] (0.41060,1.28559) circle (0.55pt);
\fill[black] (0.56161,2.35765) circle (0.55pt);
\end{tikzpicture}
\caption{The unit-distance realisation $P$ of $G_{43}$: the $43$ points and the $83$ unit edges, drawn to scale. Vertices are coloured by the $4$-colouring of Table~\ref{tab:coords}.}
\end{figure}

\begin{figure}[t]
\centering
\resizebox{0.82\textwidth}{!}{\begin{tikzpicture}[scale=1,every node/.style={circle,draw,inner sep=0pt,minimum size=4.2mm,font=\tiny}]
\node[fill=white] (v0) at (11.789,6.267) {0};
\node[fill=red!25] (v1) at (8.407,0.232) {1};
\node[fill=red!25] (v2) at (9.783,8.474) {2};
\node[fill=blue!20] (v3) at (10.483,7.491) {3};
\node[fill=blue!20] (v4) at (8.599,13.000) {4};
\node[fill=white] (v5) at (6.303,12.550) {5};
\node[fill=red!25] (v6) at (8.411,10.467) {6};
\node[fill=red!25] (v7) at (10.457,10.328) {7};
\node[fill=red!25] (v8) at (7.937,11.656) {8};
\node[fill=green!30] (v9) at (1.559,5.104) {9};
\node[fill=blue!20] (v10) at (5.026,8.360) {10};
\node[fill=white] (v11) at (7.784,6.663) {11};
\node[fill=green!30] (v12) at (5.730,10.734) {12};
\node[fill=white] (v13) at (11.356,9.724) {13};
\node[fill=blue!20] (v14) at (12.890,8.271) {14};
\node[fill=red!25] (v15) at (6.575,9.000) {15};
\node[fill=green!30] (v16) at (4.463,4.070) {16};
\node[fill=red!25] (v17) at (0.000,4.274) {17};
\node[fill=white] (v18) at (10.559,11.810) {18};
\node[fill=white] (v19) at (1.639,7.452) {19};
\node[fill=blue!20] (v20) at (5.935,3.596) {20};
\node[fill=white] (v21) at (2.819,2.578) {21};
\node[fill=red!25] (v22) at (2.610,4.029) {22};
\node[fill=white] (v23) at (2.579,11.264) {23};
\node[fill=white] (v24) at (10.358,4.247) {24};
\node[fill=white] (v25) at (5.031,0.212) {25};
\node[fill=blue!20] (v26) at (6.982,0.000) {26};
\node[fill=red!25] (v27) at (1.222,9.052) {27};
\node[fill=green!30] (v28) at (6.511,5.643) {28};
\node[fill=red!25] (v29) at (12.509,4.401) {29};
\node[fill=white] (v30) at (10.866,0.854) {30};
\node[fill=white] (v31) at (7.834,4.969) {31};
\node[fill=blue!20] (v32) at (4.465,7.512) {32};
\node[fill=blue!20] (v33) at (5.866,6.629) {33};
\node[fill=red!25] (v34) at (0.625,9.915) {34};
\node[fill=red!25] (v35) at (1.334,1.122) {35};
\node[fill=white] (v36) at (0.646,7.066) {36};
\node[fill=white] (v37) at (3.458,0.241) {37};
\node[fill=red!25] (v38) at (9.728,1.122) {38};
\node[fill=red!25] (v39) at (9.464,5.291) {39};
\node[fill=red!25] (v40) at (6.873,2.378) {40};
\node[fill=white] (v41) at (3.222,11.797) {41};
\node[fill=white] (v42) at (11.559,2.493) {42};
\begin{scope}[on background layer]
\draw[gray] (v0)--(v2);
\draw[gray] (v0)--(v28);
\draw[gray] (v0)--(v29);
\draw[gray] (v1)--(v25);
\draw[gray] (v1)--(v30);
\draw[gray] (v1)--(v31);
\draw[gray] (v2)--(v12);
\draw[gray] (v2)--(v31);
\draw[gray] (v3)--(v6);
\draw[gray] (v3)--(v24);
\draw[gray] (v3)--(v28);
\draw[gray] (v4)--(v5);
\draw[gray] (v4)--(v13);
\draw[gray] (v4)--(v15);
\draw[gray] (v5)--(v7);
\draw[gray] (v5)--(v10);
\draw[gray] (v6)--(v11);
\draw[gray] (v6)--(v12);
\draw[gray] (v7)--(v11);
\draw[gray] (v7)--(v14);
\draw[gray] (v8)--(v10);
\draw[gray] (v8)--(v13);
\draw[gray] (v8)--(v18);
\draw[gray] (v9)--(v17);
\draw[gray] (v9)--(v20);
\draw[gray] (v9)--(v21);
\draw[gray] (v9)--(v27);
\draw[gray] (v9)--(v36);
\draw[gray] (v10)--(v16);
\draw[gray] (v10)--(v23);
\draw[gray] (v10)--(v27);
\draw[gray] (v10)--(v28);
\draw[gray] (v10)--(v31);
\draw[gray] (v11)--(v16);
\draw[gray] (v11)--(v20);
\draw[gray] (v11)--(v29);
\draw[gray] (v11)--(v32);
\draw[gray] (v12)--(v23);
\draw[gray] (v12)--(v33);
\draw[gray] (v12)--(v41);
\draw[gray] (v13)--(v39);
\draw[gray] (v14)--(v18);
\draw[gray] (v14)--(v39);
\draw[gray] (v15)--(v18);
\draw[gray] (v15)--(v19);
\draw[gray] (v15)--(v31);
\draw[gray] (v15)--(v32);
\draw[gray] (v15)--(v33);
\draw[gray] (v16)--(v17);
\draw[gray] (v16)--(v25);
\draw[gray] (v16)--(v26);
\draw[gray] (v16)--(v33);
\draw[gray] (v16)--(v35);
\draw[gray] (v17)--(v19);
\draw[gray] (v19)--(v22);
\draw[gray] (v20)--(v31);
\draw[gray] (v20)--(v37);
\draw[gray] (v20)--(v39);
\draw[gray] (v21)--(v28);
\draw[gray] (v21)--(v35);
\draw[gray] (v22)--(v28);
\draw[gray] (v22)--(v37);
\draw[gray] (v23)--(v34);
\draw[gray] (v24)--(v29);
\draw[gray] (v24)--(v33);
\draw[gray] (v25)--(v40);
\draw[gray] (v26)--(v30);
\draw[gray] (v26)--(v38);
\draw[gray] (v27)--(v41);
\draw[gray] (v28)--(v32);
\draw[gray] (v28)--(v39);
\draw[gray] (v30)--(v39);
\draw[gray] (v31)--(v38);
\draw[gray] (v32)--(v34);
\draw[gray] (v32)--(v40);
\draw[gray] (v32)--(v41);
\draw[gray] (v33)--(v36);
\draw[gray] (v33)--(v39);
\draw[gray] (v34)--(v36);
\draw[gray] (v35)--(v37);
\draw[gray] (v38)--(v42);
\draw[gray] (v39)--(v42);
\draw[gray] (v40)--(v42);
\end{scope}
\end{tikzpicture}}
\caption{An abstract drawing of $G_{43}$ with vertex labels (edges are \emph{not} drawn at unit length). Vertex colours show the $4$-colouring of Table~\ref{tab:coords}; white is colour~$0$.}
\label{fig:abstract}
\end{figure}
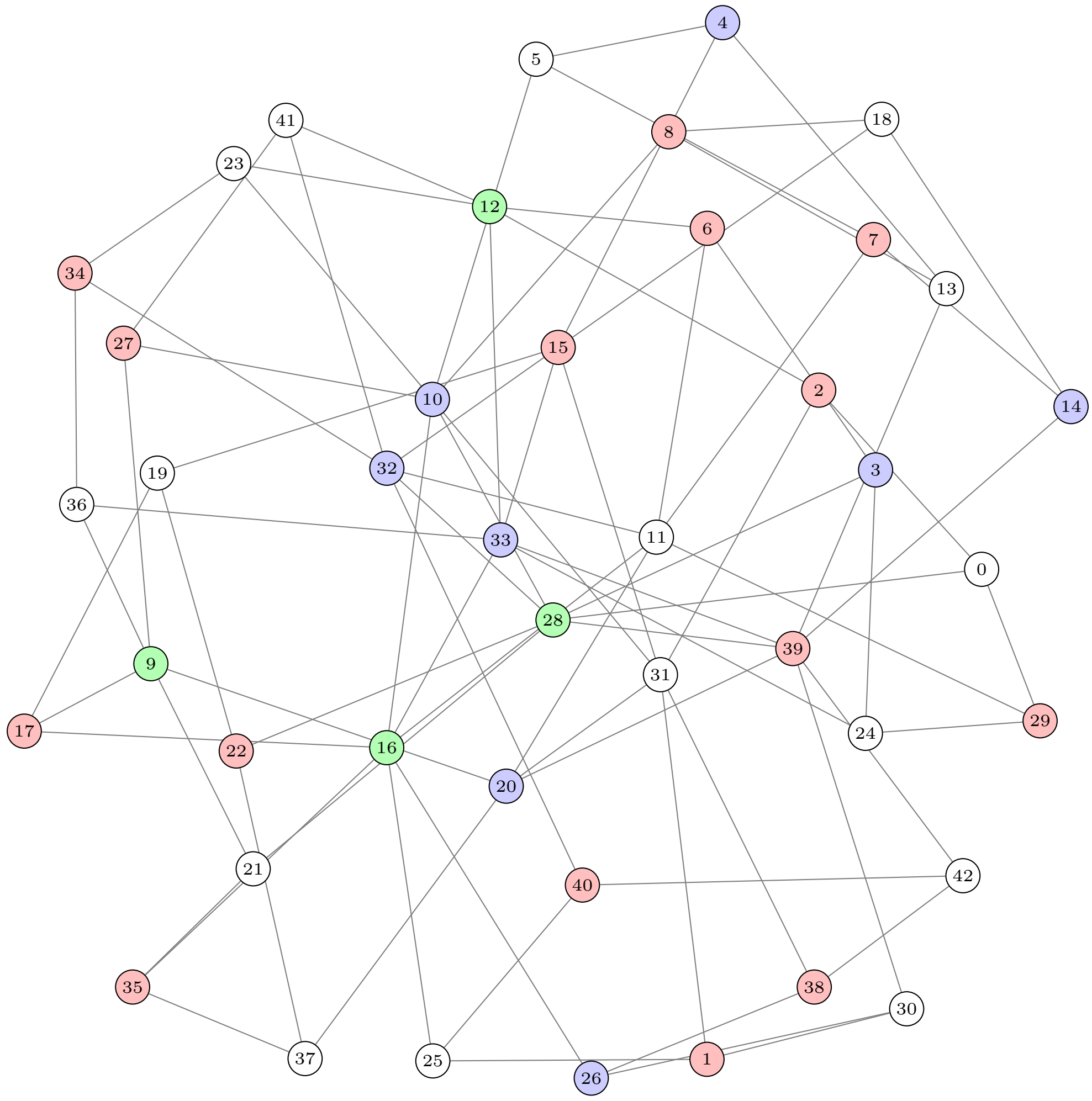

\paragraph{Exact coordinates.}
Put $p_0=(0,0)$ and $p_2=(1,0)$; the edge $\{0,2\}$ fixes position and orientation. The remaining $41$ points give $82$ unknown coordinates $t=(t_1,\dots,t_{82})$. Each of the $82$ edges other than $\{0,2\}$ gives one equation
\[
F_{uv}(t)\;=\;\|p_u-p_v\|^2-1\;=\;0 ,
\]
so $F=(F_{uv})$ is a square system of $82$ quadratic polynomials with integer coefficients. Let $c\in\mathbb{Q}^{82}$ be the vector of rational centres listed in Table~\ref{tab:coords} (given there to $16$--$17$ significant digits; the exact rational values are in the supplementary files), let $r=10^{-9}$, and let
\[
X=\{t\in\mathbb{R}^{82}:\ \|t-c\|_\infty\le r\}.
\]

\begin{lemma}\label{lem:root}
$F$ has exactly one zero $t^*$ in $X$. At $t^*$ the Jacobian $DF(t^*)$ is nonsingular. Consequently every coordinate of $t^*$ is a real algebraic number.
\end{lemma}

\begin{proof}
Let $B$ be the explicit rational $82\times 82$ matrix in the supplementary files (an approximate inverse of $DF(c)$), and set $T(t)=t-BF(t)$. Since $F$ is quadratic, each entry of $DF$ has the form $\pm2(x_u-x_v)$, and it changes by at most $4r$ on $X$; each row of $DF$ has at most four nonzero entries. Exact rational arithmetic gives
\[
\eta:=\|BF(c)\|_\infty<4.95\cdot10^{-15},\qquad
q:=\|I-B\,DF(c)\|_\infty+16r\|B\|_\infty<7.25\cdot 10^{-7}.
\]
By the mean value theorem, $\|T(s)-T(t)\|_\infty\le q\|s-t\|_\infty$ for $s,t\in X$. Also $\|T(t)-c\|_\infty\le\eta+qr<r$, so $T$ maps $X$ into itself. By the contraction mapping theorem $T$ has a unique fixed point $t^*\in X$.

Since $\|I-B\,DF(c)\|_\infty<1$, the matrix $B$ is invertible, so the fixed points of $T$ are exactly the zeros of $F$ in $X$. The same estimate holds with $DF(t^*)$ in place of $DF(c)$, so $DF(t^*)$ is nonsingular.

A nonsingular zero is an isolated point of the complex variety $\{F=0\}$. Because $F$ has rational coefficients, isolated points of this variety have algebraic coordinates.
\end{proof}

Let $P=\{p_0,\dots,p_{42}\}$ be the points given by $t^*$.

\begin{lemma}\label{lem:udg}
The points of $P$ are pairwise distinct, and the unit-distance graph of $P$ is exactly $G_{43}$.
\end{lemma}

\begin{proof}
Every edge has length $1$ by Lemma~\ref{lem:root}. For a pair $u\ne v$, interval bounds on $\|p_u-p_v\|^2$ over the box $X$, computed exactly in rational arithmetic, give
\[
\|p_u-p_v\|^2\ge 0.00173 \quad\text{for all } u\ne v,
\qquad
\bigl|\|p_u-p_v\|^2-1\bigr|\ge 0.00899\quad\text{for all non-edges } \{u,v\}.
\]
So the points are distinct and no non-edge has length $1$.
\end{proof}

\section{Proof of Theorem~\ref{thm:main}}

\emph{(1) Girth.} No edge has both ends with a common neighbour, so there are no triangles. No two vertices have two common neighbours, so there are no $4$-cycles. The cycle $0,2,31,10,28$ has length $5$.

\emph{(2) Chromatic number.} Table~\ref{tab:coords} gives a proper $4$-colouring. To show $\chi(G_{43})>3$ we encode $3$-colourability as a CNF formula with $129$ variables: one clause per vertex saying that it gets some colour, and three clauses per edge. The SAT solver \texttt{kissat} (version 4.0.4)~\cite{kissat} reports this formula unsatisfiable and produces a DRAT proof of $1916$ bytes. This proof was verified with \texttt{drat-trim}~\cite{drattrim}. No symmetry breaking is used in this formula. Independently, an exhaustive DSATUR-style backtracking search, with the ends of the edge $\{0,2\}$ fixed to colours $0$ and $1$, visits $1957$ nodes and finds no $3$-colouring.

\emph{(3) Criticality.} For each of the $43$ vertices and each of the $83$ edges, the deleted graph was $3$-coloured explicitly by backtracking. All $126$ colourings are included in the supplementary files, and each was checked to be proper.

\emph{(4) Degrees and automorphisms.} The degree sequence is read off the edge list. A complete enumeration of the automorphisms (VF2 algorithm, as implemented in \texttt{networkx}) finds only the identity.

\emph{(5) Rigidity.} The graph has $83=2\cdot 43-3$ edges. The $(2,3)$-pebble game~\cite{JacobsHendrickson} confirms that every set of $k\ge2$ vertices spans at most $2k-3$ edges, so $G_{43}$ is a Laman graph~\cite{Laman}. The rigidity matrix of the framework $(G_{43},P)$, with the edge $\{0,2\}$ pinned, is $\tfrac12 DF(t^*)$. This matrix is nonsingular by Lemma~\ref{lem:root}, so the realisation is infinitesimally rigid.
\qed

\begin{remark}
Lemma~\ref{lem:root} was checked twice, with independent code. The first check is the rational contraction certificate above, using exact \texttt{Fraction} arithmetic. The second is an interval Krawczyk test~\cite{Krawczyk,Moore} at $50$ digits, carried out with an approximate inverse computed afresh from the Newton root. That test gives $K(X)\subset\operatorname{int}X$ with $\|K(X)-c\|_\infty\le 5\cdot10^{-6}\,r$. A Newton iteration at $80$ digits from $c$ converges to a root with residual below $10^{-80}$ at distance $4.95\cdot10^{-15}$ from $c$.
\end{remark}

\begin{remark}
The realisation is tight in places: the smallest distance between two points of $P$ is about $0.0416$, and the non-edge whose length is closest to $1$ has squared length within $0.009$ of $1$. Because the system $F=0$ is square and nonsingular, $P$ is locally the only realisation of $G_{43}$ with these edge lengths, up to congruence. We have not looked for a closed-form expression of the coordinates.
\end{remark}

\section{Concluding remarks}

The minimum order of a $4$-chromatic unit-distance graph of girth $5$ is now known to lie between $21$ and $43$. It would be interesting to narrow this gap. One step would be to decide whether any of the smallest $4$-chromatic graphs of girth $5$, which have $21$ vertices (see~\cite{EG}), is a unit-distance graph.

\section*{Data availability}
The supplementary files contain:
\begin{itemize}
\item the edge list;
\item the exact rational centres and the matrix $B$;
\item the polynomial system;
\item the CNF formula and the DRAT proof;
\item two independent verification scripts.
\end{itemize}
They are available at \url{https://doi.org/10.5281/zenodo.24036368}.

\appendix
\section{Coordinates and colouring}
Decimal values of the box centres; the exact rational centres are in the supplementary files. Each free coordinate of $P$ lies within $10^{-9}$ of the value shown. The last column is a proper $4$-colouring.

\begin{table}[ht]
\centering\small
\caption{Box centres $c$ (vertices $0$ and $2$ are exact) and a $4$-colouring.}\label{tab:coords}
\begin{tabular}{rllc}
\toprule
$v$ & $x$ & $y$ & colour\\ \midrule
0 & $0$ & $0$ & 0 \\
1 & $0.4330898991676348$ & $0.17517745454218736$ & 1 \\
2 & $1$ & $0$ & 1 \\
3 & $1.2267499031050335$ & $0.0067421023916053096$ & 2 \\
4 & $1.033578176991139$ & $-0.024575643333048874$ & 2 \\
5 & $1.399381442379761$ & $-0.95526784327405057$ & 0 \\
6 & $2.1630319145538008$ & $0.35799129461480611$ & 1 \\
7 & $1.9849313733408456$ & $-0.1446315016562639$ & 1 \\
8 & $2.2506146807519531$ & $-0.11064086413120923$ & 1 \\
9 & $1.7582248930595668$ & $1.9529100593027817$ & 3 \\
10 & $1.2612996895990445$ & $0.03515299129527779$ & 2 \\
11 & $1.1655266885216353$ & $0.42858395739585431$ & 0 \\
12 & $1.3480490350752383$ & $0.93747632993222585$ & 3 \\
13 & $1.6979910603562509$ & $0.72279007732749756$ & 0 \\
14 & $1.8746120586826571$ & $0.84926469453997244$ & 2 \\
15 & $1.8741623423572826$ & $0.51710531529084758$ & 1 \\
16 & $2.1649171705193808$ & $0.46349333281728$ & 3 \\
17 & $2.5747788927340869$ & $1.3756410010504423$ & 1 \\
18 & $2.8604986057909723$ & $0.68184989163758891$ & 0 \\
19 & $1.579500949990162$ & $1.4727070439666963$ & 0 \\
20 & $1.9983241777325442$ & $0.98216172245150257$ & 2 \\
21 & $1.5938927434972694$ & $0.96650499817412217$ & 0 \\
22 & $0.6321407592436713$ & $1.7928767942379744$ & 1 \\
23 & $2.191970762407327$ & $0.40101000467329218$ & 0 \\
24 & $1.3108125298679165$ & $1.0032025756577176$ & 0 \\
25 & $1.2203468589886777$ & $0.79180252637541171$ & 0 \\
26 & $1.2698663108079469$ & $0.90945741043993222$ & 2 \\
27 & $1.3769681292103328$ & $1.028440871055164$ & 1 \\
28 & $0.6090344347145451$ & $0.79314378099556138$ & 3 \\
29 & $0.31222841357043618$ & $0.95000706195232509$ & 1 \\
30 & $0.30353260224069506$ & $1.1667493918281235$ & 0 \\
31 & $0.99848327194346864$ & $0.99999884976733966$ & 0 \\
32 & $1.4040613435456299$ & $1.3997179431839042$ & 2 \\
33 & $2.1988796223628104$ & $1.4629164423483652$ & 2 \\
34 & $2.4038321614271041$ & $1.3783097271767739$ & 1 \\
35 & $2.4906931595836426$ & $1.4089403196701669$ & 1 \\
36 & $2.6805804484767219$ & $2.3392521737442582$ & 0 \\
37 & $1.6252619207765713$ & $1.9099680343695213$ & 0 \\
38 & $1.4474009678085034$ & $1.8935719594451255$ & 1 \\
39 & $1.2075427605754268$ & $1.5942603666903059$ & 1 \\
40 & $0.40450088711542082$ & $1.3700717977172157$ & 1 \\
41 & $0.41059561911934073$ & $1.2855870779889529$ & 0 \\
42 & $0.56160838256353762$ & $2.3576533056820117$ & 0 \\
\bottomrule
\end{tabular}
\end{table}

\end{document}